\documentclass{amsart}

\usepackage{amsfonts}
\usepackage{amsmath}
\usepackage{amsthm}
\usepackage{amssymb}
\usepackage[english]{babel}
\usepackage{graphics}
\usepackage{latexsym}
\usepackage{longtable} 
\usepackage{mathrsfs} 
\usepackage{graphicx}
\usepackage{wrapfig} 

\usepackage{hyperref}

\usepackage{enumitem}
\usepackage{esint}
\usepackage{mathcomp}
\usepackage{stmaryrd}
\usepackage{wasysym}

\newtheorem{theorem}{Theorem}
\newtheorem{corollary}[theorem]{Corollary}
\newtheorem{lemma}[theorem]{Lemma}

\newtheorem{claim}[theorem]{Claim}
\newtheorem{example}[theorem]{Example}

\theoremstyle{definition}
\newtheorem{definition}[theorem]{Definition}
\newtheorem{remark}[theorem]{Remark}

\renewcommand{\d}{\mathrm{d}}

\newcommand{\mL}{\mathcal{L}}

\newcommand{\mF}{\mathcal{F}}

\newcommand{\mX}{\mathfrak{X}}

\newcommand{\mM}{\mathcal{M}}

\newcommand{\E}{\mathrm{E}}
\newcommand{\D}{\mathrm{D}}

\newcommand{\W}{\mathcal{W}}

\newcommand{\K}{\mathrm{K}}

\newcommand{\G}{\mathrm{G}}

\newcommand{\R}{\mathbb{R}}

\newcommand{\N}{\mathbb{N}}

\newcommand{\mB}{\mathbb{B}}

\newcommand{\X}{\mathrm{X}}

\newcommand{\Q}{\mathrm{Q}}

\newcommand{\ms}{\medskip}

\newcommand{\al}{\alpha}

\newcommand{\Ga}{\Gamma}

\newcommand{\e}{\varepsilon}
\newcommand{\si}{\sigma}
\newcommand{\Si}{\Sigma}
\newcommand{\la}{\lambda}
\newcommand{\La}{\Lambda}

\newcommand{\Om}{\Omega}

\newcommand{\av}{-\hspace{-10.5pt}\displaystyle\int}

\newcommand{\weak }{\, -\!\!\!\!\!-\!\!\!\!\rightharpoonup}
\newcommand{\weakstar }{ \overset{\, *_{\phantom{|}}}{{\smash{\, -\!\!\!\!-\!\!\!\!\rightharpoonup}}\, } }

\newcommand{\larrow}{\longrightarrow}

\newcommand{\LL}{\text{\LARGE$\llcorner$}}
\newcommand{\p}{\partial}
\newcommand{\sub}{\subseteq}

\newcommand{\by}{\times}

\renewcommand{\div}{\mathrm{div}}

\renewcommand{\ker}{\mathrm{ker}}

\newcommand{\bt}{\begin{theorem}}\newcommand{\et}{\end{theorem}}
\newcommand{\bd}{\begin{definition}}\newcommand{\ed}{\end{definition}}
\newcommand{\bl}{\begin{lemma}}\newcommand{\el}{\end{lemma}}
\newcommand{\beq}{\begin{equation}}\newcommand{\eeq}{\end{equation}}
\newcommand{\bc}{\begin{claim}}\newcommand{\ec}{\end{claim}}
\newcommand{\bex}{\begin{example}}\newcommand{\eex}{\end{example}}
\newcommand{\bcor}{\begin{corollary}}\newcommand{\ecor}{\end{corollary}}
\newcommand{\bp}{\begin{proof}}\newcommand{\ep}{\end{proof}}

\newcommand{\BPL}{\medskip \noindent \textbf{Proof of Lemma} }

\newcommand{\BPT}{\medskip \noindent \textbf{Proof of Theorem} }

\numberwithin{equation}{section}

\begin{document}

\title[Variational problems in $L^\infty$ constrained by the Navier-Stokes]{Vectorial variational problems in $L^\infty$ constrained by the Navier-Stokes equations 
}
 
\author{Ed Clark}

\address{E. C., Department of Mathematics and Statistics, University of Reading, Whiteknights Campus, Pepper Lane, Reading RG6 6AX, United Kingdom}

\email{e.d.clark@pgr.reading.ac.uk}
 
\author{Nikos Katzourakis}

\address[Corresponding author]{N. K., Department of Mathematics and Statistics, University of Reading, Whiteknights Campus, Pepper Lane, Reading RG6 6AX, United Kingdom}

\email{n.katzourakis@reading.ac.uk}
 
\author{Boris Muha}

\address{B. M., Department of Mathematics, Faculty of Science, University of Zagreb, Zagreb, Trg Republike Hrvatske 14, 10000, Croatia}

\email{borism@math.hr}

\thanks{E.C.\ has been financially supported through the UK EPSRC scholarship GS19-055}
\thanks{B.M.\ has been partially financially supported through the Croatian Science Foundation project  IP-2019-04-1140}
  



\begin{abstract} We study a minimisation problem in $L^p$ and $L^\infty$ for certain cost functionals, where the class of admissible mappings is constrained by the Navier-Stokes equations. Problems of this type are motivated by variational data assimilation for atmospheric flows arising in weather forecasting. Herein we establish the existence of PDE-constrained minimisers for all $p$, and also that $L^p$ minimisers converge to $L^\infty$ minimisers as $p\to\infty$. We further show that $L^p$ minimisers solve an Euler-Lagrange system. Finally, all special $L^\infty$ minimisers constructed via  approximation by $L^p$ minimisers are shown to solve a divergence PDE system involving measure coefficients, which is a divergence-form counterpart of the corresponding non-divergence Aronsson-Euler system. 
\end{abstract}

\maketitle

\tableofcontents

\section{Introduction and main results}   \label{Section1}

Let $\Om \sub \R^n$ be an open bounded set and let also $n\geq 2$ and $\nu,T>0$. Consider the Navier-Stokes equations
\beq 
\label{1.1}
\left\{\ 
\begin{array}{lll}
 & \ \  \p_t u\,-\nu\Delta u\,+ (u\cdot \D)u\, +\D \mathrm p \,-f =\, y, & \ \ \text{ in }\Om\times (0,T), \ms
\\
 &\ \     \div \, u\, =0, & \ \ \text{ in }\Om\times (0,T), \ms
\\
  & \ \  u(\cdot,0) =u_0, & \ \ \text{ on }\Om,  \ms
\\
 & \ \  u =\, 0, & \ \ \text{ on }\p\Om\times (0,T),
\end{array}
\right.
\eeq
and for brevity let us henceforth symbolise $\nabla u := (\D u,\p_t u)$ and $\Om_T := \Om\by (0,T)$, where $\D u=(\p_{x_1} u,...,\p_{x_n} u) \in \R^{n \by n}$ symbolises the spatial gradient.  The system of PDEs \eqref{1.1} describes the velocity $u: \Om_T \larrow \R^n$ and the pressure $\mathrm p: \Om_T \larrow \R$ of a flow, for some given initial data $u_0 : \Om \larrow \R^n$ with source $f: \Om_T \larrow \R^n$. Here the map $y: \Om_T \larrow \R^n$ is a parameter and should be understood as a (deterministic) noise or error. Let also $N\in \N$ and suppose we are given a mapping
\beq \label{1.2}
\Q \ : \ \ \Om_T \by \big(\R^n \by {\R}^{(n+1 )\by n} \times \R\big) \larrow {\R}^N.
\eeq
A problem of interest in the geosciences, in particular in data assimilation for atmospheric flows in relation to weather forecasting (see e.g.\ \cite{B1,B2,B3}), can be formulated as follows: find solutions $(u, \mathrm p)$ to \eqref{1.1} such that, in an appropriate sense, 
\beq 
\label{1.3}
\left\{ \ \ 
\begin{split}
y \, \approx \, 0,& 
\\
\Q(\cdot,\cdot, u, \nabla u, \mathrm p) -q \, \approx \, 0,& \ \ \ \ \ \ \ 
\end{split}
\right.
\eeq
where $q: \Om_T \larrow \R^N$ is a vector of given measurable ``data" arising from some specific measurements, taken through the ``observation operator" $\Q$ of \eqref{1.2}. In \eqref{1.1} and \eqref{1.3}, $y$ represents an error in the measurements which forces the Navier-Stokes equations to be satisfied only approximately for solenoidal (divergence-free) vector fields. Namely, we are looking for solutions to \eqref{1.1} such that simultaneously the error $y$ vanishes, and also the measurements $q$ match the prediction of the solution $(u,\mathrm p)$ through the observation operator \eqref{1.2}. In application, $\Q$ is typically some component (e.g.\ linear projection or nonlinear submersion) of the atmospheric flow that we can observe. Unfortunately, the data fitting problem \eqref{1.3} is severely ill-posed; an exact solution may well not exist, and even if it does, it may not be unique. 

In this paper, inspired by the methodology of data assimilation, especially variational data assimilation in continuous time (for relevant works we refer e.g.\ to \cite{BOT, BKO, DPV, FLT, FMT, Ko, LP, SW}), we seek to minimise the misfit functional
\[
(u,\mathrm p, y) \, \mapsto \, (1-\lambda)\big\|\Q(\cdot,\cdot, u,\nabla u, \mathrm p )-q\big\| \, +\, \lambda\big\|y\big\|
\]
over all admissible triplets $(u,\mathrm p, y)$ which satisfy \eqref{1.1}, for a fixed weight $\lambda \in (0,1)$. The role of this weight is to obtain essentially a Pareto family of extremals, one for each value $\la$, even though in this paper we do not pursue further this viewpoint of vector-valued minimisation (the interested reader may e.g.\ consult \cite{CHY}). The standard approach to data assimilation is to use Hilbert space methods (or least squares in the discrete case), hence minimising in $L^2$. The novelty of our approach, which is also justified from the viewpoint of applications, is to consider instead \emph{minimisation in $L^\infty$}, namely by interpreting the norms above as $L^\infty$ ones (or maxima in the discrete case). There is a significant advantage of considering a min-max problem instead of minimising the squared averages: the misfit becomes uniformly small throughout the space-time domain $\Om_T$ and not just on average, hence large ``spikes" of deviations from zero are at the outset excluded.

When one passes from a variational problem for an integral norm to one for the supremum norm, even though this is justified from the viewpoint of desired outputs, it poses some serious challenges. The $L^\infty$ norm is neither differentiable nor strictly convex, and the space $L^\infty$ is neither reflexive nor separable. Additionally, with respect to the domain argument, the $L^\infty$ norm is not additive but only sub-additive. Further, one would also need estimates for \eqref{1.1} in appropriate subspaces of $L^\infty$ for weakly differentiable functions, which, to the best of our knowledge, do not exist even for linear strongly elliptic systems (see e.g.\ \cite{GM}). Even then, if one somehow solves the $L^\infty$ minimisation problem (by using, for instance, the direct method of the Calculus of Variations as in \cite{D}, under the appropriate quasiconvexity assumptions for $|\Q-q|+|y|$ as in \cite{BJW2}), the analogue of the Euler-Lagrange equations for the $L^\infty$ problem cannot be derived directly by perturbation/sensitivity methods due to the lack of smoothness of the $L^\infty$ norm.

In this paper, to overcome the difficulties described above, we follow the methodology of the relatively new field of Calculus of Variations in $L^\infty$ (see e.g.\ \cite{C, K0} for a general introduction to the scalar-valued theory), and in particular the ideas from \cite{K2, K3, K4, KM} involving higher order and vectorial problems, as well as problems involving PDE-constraints, which have only recently started being investigated. To this end, we follow the approach of solving the desired $L^\infty$ variational problem by solving respective approximating $L^p$ variational problems for all $p$, and obtain appropriate compactness estimates which allow to pass to the limit as $p\to \infty$. The case of finite $p>2$ studied herein is also of independent interest, especially for numerical discretisation schemes in $L^\infty$ (see e.g.\ \cite{KP1, KP2}), but in this paper we treat it mostly as an approximation device to solve efficiently the $L^\infty$ problem. The idea of this approach is based on the observation that, for a fixed essentially bounded function on a domain of finite measure, the $L^p$ norm tends to the $L^\infty$ norm of the function as $p\to\infty$. 

In order to state our hypotheses and main results, let us set
\beq \label{1.4}
\K \big(x,t,\eta,A,a,r\big) \,:=\, \Q\big(x,t,\eta,A,a,r\big)\, -\, q(x,t)
\eeq
(note that in \eqref{1.4} $(x,t)\in\Om_T$ is treated as two arguments and the two arguments $(A,a)$ are for $\nabla u =(\D u,\p_t u)$, which we conveniently display abbreviated as one) and, by considering the (strong) divergence operator $\div : W^{1,1}(\Om;\R^n) \larrow L^1(\Om)$, we henceforth assume that
\beq 
\label{1.5}
\left\{ \ \ 
\begin{array}{ll}
(a) &\Om \text{ is bounded and has $C^2$ boundary }\p\Om, \ms\\
(b) &u_0\in \big(W^{2,{\infty}}\cap W^{1,\infty}_0\big)({\Om};{\R}^n)\cap \ker(\div),  \ms\\
(c) & f\in L^{\infty}({\Om}_T;{\R}^n) \ \ \ \& \ \ \  q\in L^{\infty}({\Om}_T;{\R}^N), \ms\\
(d) &\K(x,t,\cdot,\cdot,\cdot,\cdot) \text{ is $C^1$ for almost every }(x,t), \ms\\
(e) &\K(\cdot,\cdot,\eta,A,a,r) \text{ is $L^\infty$ for all }(\eta,A,a,r), \ms\\
(f) & |\K(x,t,\eta,A,\cdot,\cdot)|^2 \text{ is convex for all }(x,t,\eta,A).
\end{array}
\right.
\eeq
Then, for any $p\in (1,\infty)$, we define the $L^p$ misfit $\E_p : \mathfrak{X}^{p}({\Om}_T) \larrow \R$ by setting
\beq 
\label{1.6}
{\E}_p\big(u, \mathrm p  ,y \big) \,:= \, 
(1-\lambda){\big\|\K(\cdot, u,\nabla u, \mathrm p )\big\|}_{\dot{L}^p({\Om}_T)}
 + \, \lambda \| y \|_{\dot{L}^p({\Om}_T)}.
\eeq
We note that in \eqref{1.6} and subsequently, the dotted $\dot{L}^p$ quantities are regularisations of the respective norms at the origin, obtained by regularising the Euclidean norm in the respective target space:
\beq
\label{1.7}
\| h \|_{\dot{L}^p({\Om}_T)} :=\, \big\| |h|_{(p)} \big\|_{L^p({\Om}_T)}, \ \ \ \   | \cdot |_{(p)} \,:=\, \sqrt{| \cdot |^2 + p^{-2}}.
\eeq
Further, since we will only be dealing with finite measures, we will always be using the normalised $L^p$ norms in which we replace the integral over the domain with the respective average, for example for $L^p(\Om_T)$ with the $(n+1)$-Lebesgue measure, the norm will be
\[
\|h \|_{L^p(\Om_T)} := \left(\, \av_{\Om_T}|h|^p\, \mathrm d \mL^{n+1}\! \right)^{\!\!1/p}.
\]
The admissible minimisation class $\mathfrak{X}^{p}({\Om}_T)$ over which $\E_p$ is considered, is defined as follows:
\beq
\label{1.8}
\mathfrak{X}^p(\Om_T) \,:=\, \Big\{ (u,\mathrm p,y) \in \W^p(\Om_T) \ : \ (u,\mathrm p,y)\text{ satisfies weakly }\eqref{1.1}\Big\},
\eeq
where
\beq
\label{1.9}
\W^p(\Om_T) \, :=\, W^{2,1;p}_{\mathrm L,\si}(\Om_T;\R^n) \by W^{1,0;p}_\sharp(\Om_T) \by L^p(\Om_T;\R^n).
\eeq
The rather complicated functional spaces appearing in \eqref{1.9} are defined as follows. The space $W^{2,1;p}_{\mathrm L,\si}(\Om_T;\R^n)$ consists of solenoidal maps which are $W^{2,p}$ in space and $W^{1,p}$ in time, and also laterally vanishing on $\p \Om \by (0,T)$:
\beq
\label{1.10}
\left\{ 
\begin{split}
W^{2,1;p}_{\mathrm L,\si}(\Om_T;\R^n)\, &:=\, L^p\big((0,T);W^{2,p}_{0,\si}(\Om;\R^n)\big) \bigcap W^{1,p}\big((0,T);L^p(\Om;\R^n)\big),
\\
W^{2,p}_{0,\si}(\Om;\R^n) \, &:=\, \big( W^{2,p} \cap W^{1,p}_0\big) (\Om;\R^n) \cap \ker(\div).  
\end{split}
\right. \!\!
\eeq
The space $W^{1,0;p}_\sharp (\Om_T)$ consists of scalar-valued functions which are $W^{1,p}$ in space with zero average, and $L^p$ in time:
\beq
\label{1.11}
\left\{ \ \ 
\begin{split}
W^{1,0;p}_\sharp (\Om_T)\, &:=\, L^p\big((0,T);W^{1,p}_\sharp(\Om)\big),
\\
W^{1,p}_\sharp(\Om)\,& :=\, \left\{ g \in W^{1,p}(\Om)\ : \ \int_\Om g\,\d\mL^n =0 \right\}.
\end{split}
\right.
\eeq
The associated norms in these spaces are the expected ones, namely
\beq
\label{1.12}
\left\{
\ \ 
\begin{split}
\|v\|_{W^{2,1;p}_{\mathrm L,\si}(\Om_T)} \, &:=\, \|v\|_{L^p(\Om_T)} \, +\, \|\nabla v\|_{L^p(\Om_T)} \,+\, \|\D^2v\|_{L^p(\Om_T)},
\\
\| g \|_{W^{1,0;p}_{\sharp}(\Om_T)} \, &:=\, \|g\|_{L^p(\Om_T)} \, +\, \|\D g \|_{L^p(\Om_T)}.
\end{split}
\right.
\eeq
Note also that the divergence-free condition for $u$ in \eqref{1.1} has now been incorporated in the functional space $W^{2,1;p}_{\mathrm L,\si}(\Om_T)$. Finally, the $L^\infty$ misfit $\E_\infty : \mathfrak X^\infty({\Om}_T) \larrow \R$ is defined by setting
\beq 
\label{1.13}
\E_\infty \big(u, \mathrm p  ,y \big) \,:= \, 
(1-\lambda){\big\|\K(\cdot,\cdot, u,\nabla u, \mathrm p )\big\|}_{L^\infty({\Om}_T)}
 + \, \lambda \| y \|_{L^\infty({\Om}_T)},
\eeq
where the admissible class $\mathfrak X^\infty(\Om_T)$ is given by
\beq
\label{1.14}
\mathfrak X^\infty(\Om_T)\, :=\, \bigcap_{1<p<\infty} \mX^p(\Om_T).
\eeq
Note that the natural topology of $\mX^\infty(\Om_T)$ is not induced by a complete norm in a Banach space, but instead its topology is defined as the locally convex topology induced from the ambient Frech\'et space $\smash{\bigcap_{1<p<\infty}\W^p(\Om_T)}$. Notwithstanding, no difficulties will arise from this fact, which is a necessity that stems from the lack of $W^{2,\infty}$-$W^{1,\infty}$ estimates for \eqref{1.1}. In particular, $\mathfrak X^\infty(\Om_T)$ is \emph{not} obtained by considering the strictly smaller Cartesian product space
\[
\W^\infty(\Om_T) \, =\, W^{2,1;\infty}_{\mathrm L,\si}(\Om_T) \by \smash{W^{1,0;\infty}_{\sharp}(\Om_T)} \by L^\infty(\Om_T;\R^n). 
\]
However, we will assume that the solution $(u,\mathrm p)$ to \eqref{1.1} is strong and satisfies $W^{2,p}$-$W^{1,p}$ estimates for any finite $p$. This is deducible under assumption \eqref{1.5} in the case of $n=2$ (see e.g.\ \cite{Ge,S}), and also under smallness conditions on $u_0$ in any dimension $n\geq 3$ (see e.g.\ \cite{Amann2,Gi,So}). Hence, our additional hypothesis is
\beq
\label{1.15}
\left\{ \ \ 
\begin{split}
& \text{Either}
\\
& \hspace{130pt} n=2, \phantom{_\big]}
\\
&\text{or, $n\geq3$ but for any $p\in(1,\infty)$, exists $C>0$ depending only on $p$ and}
\\
&\text{on $\p\Om$, $T$, $\| u_0\|_{L^2(\Om)},\| f\|_{L^2(\Om_T)}$, such that}
\\
&\ \ \|u\|_{W^{2,1;p}_{\mathrm L,\si}(\Om_T)} +\, \|\mathrm p \|_{W^{1,0;p}_{\sharp}(\Om_T)} \, \leq \, C\big(\| u_0\|_{W^{2-\frac{2}{p},p}(\Om)}  +\, \| f\|_{L^p(\Om_T)}\big), \phantom{\bigg]}
\\
& \text{when $(u,\mathrm p)$ solves weakly \eqref{1.1} with $y\equiv 0$.} 
\end{split}
\right. \!\!\!
\eeq
Assumption \eqref{1.15}, albeit restrictive, is compatible with situations of interest in weather forecasting (see e.g.\ \cite{B1,B2,B3}). Our first main result concerns the existence of $\E_p$-minimisers in $\mathfrak X^p(\Om_T)$, the existence of $\E_\infty$-minimisers in $\mathfrak X^\infty(\Om_T)$ and the approximability of the latter by the former as $p\to \infty$. 

\begin{theorem}[$\E_\infty$-minimisers, $\E_p$-minimisers \& convergence as $p\to\infty$] \label{theorem1} 
Suppose that \eqref{1.5} and \eqref{1.15} hold true. Then, for any $p \in (n+2,\infty]$, the functional $\E_p$ (given by \eqref{1.6} for $p<\infty$ and by \eqref{1.13} for $p=\infty$) has a constrained minimiser $(u _p,\mathrm p_p ,y_p)$ in the admissible class $\mathfrak X^p(\Om_T)$:
\beq \label{1.16}
\E_p\big(u_p, \mathrm p_p, y_p \big)\, =\, \inf\Big\{\E_p\big(u, \mathrm p, y \big)\, : \ \big(u, \mathrm p, y \big) \in \mX^p(\Om_T) \Big\}.
\eeq
Additionally, there exists a subsequence of indices $(p_j)_1^\infty$ such that the sequence of respective $\E_{p_j}$-minimisers $(u _{p_j},\mathrm p_{p_j} ,y_{p_j})$ satisfies $(u_p, \mathrm p_p, y_p) \weak (u_\infty, \mathrm p_\infty, y_\infty)$ in $\W^q(\Om_T)$ for any $q\in (1,\infty)$, as $p_j\to\infty$. Additionally,
\beq
\label{1.17}
\left\{ \ \
\begin{array}{ll}
u_{p} \weak u_\infty, &  \ \ \text{in } W^{2,1;q}_{\mathrm L,\si}(\Om_T;\R^n), \smallskip
\\
u_{p} \larrow u_\infty, &  \ \ \text{in } C \big(\overline{\Om_T};\R^n \big), \smallskip
\\
\D u_{p} \larrow \D u_\infty, &  \ \ \text{in } C \big(\overline{\Om_T};\R^{n\by n} \big), \smallskip
\\
\mathrm  p_p \weak \mathrm p_\infty, &  \ \ \text{in } W^{1,0;q}_{\#}(\Om_T;\R^{n}),  \smallskip
\\
y_p \weak y_\infty, & \ \  \text{in }L^q(\Om_T),
\end{array}
\right.
\eeq
for any $q\in (1,\infty)$, and also
\beq
\label{1.18}
\E_p (u_{p}, \mathrm p_{p} ,y_{p} ) \larrow  \E_\infty (u_\infty, \mathrm p_\infty ,y_\infty )
\eeq
as $p_j\to\infty$.
\end{theorem}
Given the existence of constrained minimisers established by Theorem \ref{theorem1} above, the next natural question concerns the existence of necessary conditions in the form of PDEs governing the constrained minimisers. We first consider the case of $p<\infty$. Unsurprisingly, the PDE constraint of \eqref{1.1} used in defining \eqref{1.8} gives rise to a generalised Lagrange multiplier in the Euler-Lagrange equations, obtained by utilising well-known results on the Kuhn-Tucker theory from \cite{Z}. Interestingly, however, the incorporation of the solenoidality constraint into the functional space (recall \eqref{1.10}), allows us to have only one generalised multiplier corresponding only to the parabolic system in \eqref{1.1}, instead of two.

To state our second main result, we first need to introduce some notation. For any $M\in\N$ and $p\in (1,\infty)$, we define the operator 
\[
\mathfrak{M}_p \ : \ \ L^p(\Om_T;\R^M) \larrow L^{p'}(\Om_T;\R^M), 
\]
where $p':=p/(p-1)$, by setting
\beq 
\label{1.19}
\mathfrak{M}_p(V)\, :=\, \frac{ |V|_{(p)}^{p-2}\, V}{  \big(\|V \|_{\dot{L}^p(\Om_T)}\big)^{p-1}} .
\eeq
Here $|\cdot |_{(p)}$ is the regularisation of the Euclidean norm of $\R^M$, as defined in \eqref{1.7}. By H\"older's inequality it is immediate to verify that (for the normalised $\smash{L^{p'}}$ norm) we actually have
\[
\big\| \mathfrak{M}_p(V) \big\|_{L^{p'}(\Om_T)} \, \leq \, 1, 
\]
and therefore $\mathfrak{M}_p$ is valued in the unit ball of $\smash{L^{p'}(\Om_T;\R^M)}$. Further, for brevity we will use the notation
\beq
 \label{K[u,p]}
 \left\{ \ \ 
 \begin{split}
\K[u,\mathrm p]\,:= &\, \K\big(\cdot,\cdot,u, \nabla u, \mathrm p\big),
\\
\K_\eta[u,\mathrm p]\,:= & \, \K_\eta\big(\cdot, \cdot,u, \nabla u, \mathrm p\big),
\\
\K_{(A,a)}[u,\mathrm p]\,:= &\, \K_{(A,a)}\big(\cdot, \cdot,u, \nabla u, \mathrm p\big),
\\
\K_r[u,\mathrm p]\,:= &\, \K_r\big(\cdot,\cdot, u, \nabla u, \mathrm p\big),
\end{split}
\right.
\eeq
for $\K$ and its partial derivatives $\K_\eta,\K_{(A,a)},\K_r$ with respect to the arguments for $u,\nabla u$ and $\mathrm p$ respectively.

\begin{theorem}[Variational Equations in $L^p$] \label{theorem2} 
Suppose that \eqref{1.5} and \eqref{1.15} hold true. Then, for any $p\in(n+2,\infty)$, there exists a Lagrange multiplier
\beq
\label{1.24}
\Psi_p \, \in \, \Big(W^{2-\frac{2}{p},p}_{0,\si} (\Om;\R^{n})\Big)^*
\eeq
associated with the constrained minimisation problem \eqref{1.16}, such that the minimising triplet $(u_{p}, \mathrm p_{p} ,y_{p} ) \in \mX^p(\Om_T)$ satisfies the relations
\beq
\label{1.22}
\left\{ \ 
\begin{split}
 &(1-\lambda)\int_{\Om_T} \Big(\K_\eta [u_p,\mathrm p_p]\cdot u \, +\, \K_{(A,a)}[u_p,\mathrm p_p] : \nabla u \Big)\cdot \mathfrak{M}_p \big( \K[u_p,\mathrm p_p] \big) \, \d\mL^{n+1}
 \\
 &= 
 -\la  \int _{\Om_T}\Big(\p_t u - \nu\Delta u + (u\cdot \D)u_p+( u_p\cdot \D)u \Big)\cdot \mathfrak{M}_p(y_p) \, \d {\mathcal{L}}^{n+1}\, + \, \big\langle \Psi_p , u(\cdot,0) \big\rangle,
\end{split}
\right.
\eeq
\beq
\label{1.23}
\begin{split}
&(1-\lambda)\int_{\Om_T}   \K_r[u_p,\mathrm p_p] \, \mathrm p \cdot \mathfrak{M}_p \big( \K[u_p,\mathrm p_p] \big)  \, \d \mathcal{L}^{n+1}\, =  -\la \int_{\Om_T}\D \mathrm p \cdot \mathfrak{M}_p(y_p) \, \d \mL^{n+1}
\\
\end{split}
\eeq
for all test mappings
\[
(u, \mathrm p) \, \in \, W^{2,1;p}_{\mathrm L,\si}(\Om_T;\R^n) \by W^{1,0;p}_\sharp (\Om_T),
\] 
where the operators $\K,\K_\eta,\K_{(A,a)},\K_r$ are given by \eqref{K[u,p]}.
\end{theorem}

Now we consider the case of $p=\infty$. For this extreme case, which is obtained by an appropriate passage to limits as $p\to\infty$ in Theorem \ref{theorem2}, we need to assume additionally that the operator $\K[u,\mathrm p ]$ does not depend on $(\p_t u, \mathrm p)$, hence in this case we will symbolise
\beq
 \label{K[u]}
 \left\{ \ \ 
 \begin{split}
\K[u]\,:= &\, \K\big(\cdot,\cdot,u, \D u \big),
\\
\K_\eta[u ]\,:= & \, \K_\eta\big(\cdot, \cdot,u, \D u\big),
\\
\K_{A}[u]\,:= &\, \K_A \big(\cdot, \cdot,u, \D u \big),
\end{split}
\right.
\eeq
for $\K$ and its partial derivatives $\K_\eta,\K_{A}$ with respect to the arguments for $u,\D u$ respectively, all of which will also need to be assumed to be continuous. We note that, when $p=\infty$, there is no direct analogue of the divergence structure Euler Lagrange equations. Instead, one of the central points of Calculus of Variations in $L^\infty$ is that Aronsson-Euler PDE systems may be derived, under appropriate (stringent) assumptions. Even in the unconstrained case, these PDE systems are always non-divergence and even fully nonlinear and with discontinuous coefficients (see e.g.\ \cite{AK1, AK2, CKP, K1, KP2}). The case of $L^\infty$ problems involving only first order derivative of scalar-valued functions is nowadays a well established field which originated from the work of Aronsson in the 1960 \cite{A1,A2}, today largely interconnected to the theory of Viscosity Solution to nonlinear elliptic PDE (for a general pedagogical introduction see e.g.\ \cite{C,K0}). However, vectorial and higher $L^\infty$ variational problems  involving constraints, have only recently began being explored (see \cite{K3,K4}, but also the relevant earlier contributions \cite{AP,AB,BJ}). For several interesting developments on $L^\infty$ variational problems we refer the interested reader to \cite{BBJ, BJW1, BN, BP, CDP, DPV, GNP, KZ, MWZ, PZ, RZ}.

In this paper, motivated by recent progress on higher order and on constrained $L^\infty$ variational problems made in \cite{KM} by the author jointly with Moser and by the author in \cite{K3,K4} (inspired by earlier contributions by Moser and Schwetlick deployed in a geometric setting in \cite{MS}), we follow a slightly different approach which does not lead an Aronsson-Euler type system; instead, it leads to a \emph{divergence structure} PDE system. However, there is a toll to be paid, as the divergence PDEs arising as necessary conditions involve measures as auxiliary parameters whose determination becomes part of the problem. Notwithstanding, the central point of this idea is to use a scaling in the Euler-Lagrange equations before letting $p\to\infty$, which is different from the scaling used to (formally) derive the Aronsson-Euler equations as $p\to\infty$. 

In the light of the above comments, our final main result concerns the satisfaction of necessary PDE conditions for the PDE-constrained minimisers in $L^\infty$ constructed in Theorem \ref{theorem1}, and reads as follows.

\begin{theorem}[Variational Equations in $L^\infty$] \label{theorem3} 
Suppose that \eqref{1.5} and \eqref{1.15} hold true, and that additionally $\K$ does not depend on $(\p_tu,\mathrm p)$ with $\K,\K_\eta,\K_A$ in \eqref{K[u]} being continuous on $\overline{\Om_T}\by \R^n \by \R^{n\by n}$. Then, there exists a linear functional
\beq
\label{1.21}
\Psi_\infty  \in   \bigcap_{r>n+2} \Big(W^{2-\frac{2}{r},r}_{0,\si} (\Om;\R^n)\Big)^* 
\eeq
which is a Lagrange multiplier associated with the constrained minimisation problem \eqref{1.16} for $p=\infty$. There also exist vector measures
\beq
\label{1.25}
\Sigma_\infty \in \mM\big(\overline{\Om_T};\R^N \big), \ \ \ \sigma_\infty \in \mM\big(\overline{\Om_T};\R^n \big)
\eeq
such that the minimising triplet $(u_\infty, \mathrm p_\infty ,y_\infty) \in \mX^\infty(\Om_T)$ satisfies the relations
\beq
\label{1.26}
\left\{ \ 
\begin{split}
 &(1-\lambda)\int_{\overline{\Om_T}} \Big(\K_\eta [u_\infty]\cdot u \, +\, \K_A[u_\infty] : \D u \Big) \cdot \d\Sigma_\infty
 \\
 &= 
 -\la  \int _{\overline{\Om_T}}\Big(\p_t u - \nu\Delta u + (u\cdot \D)u_\infty+( u_\infty\cdot \D)u \Big)\cdot \d \sigma_\infty\, + \,\big\langle \Psi_\infty , u(\cdot,0) \big\rangle,
\end{split}
\right.
\eeq
\beq
\label{1.27}
\int_{\overline{\Om_T}}\, \D \mathrm p \cdot \d \sigma_\infty \, =\, 0,
\eeq
for all test mappings
\[
(u, \mathrm p) \, \in \, \Big(W^{2,1;\infty}_{\mathrm L,\si}(\Om_T;\R^n) \cap C^2\big(\overline{\Om_T};\R^n \big) \Big) \by \Big( W^{1,0;\infty}_\sharp (\Om_T) \cap C^1\big(\overline{\Om_T}\big) \Big).
\] 
Further, the multiplier $\Psi_\infty$ and the measures $\Sigma_\infty,\sigma_\infty$ can be approximated as follows:
\beq
\label{1.28}
\left\{ \ \
\begin{split}
&\Psi_p \, \weakstar  \, \Psi_\infty,  \ \, \text{ in }\smash{\big( W^{2-{2}/{r},r}_{0,\si} (\Om;\R^n) \big)^*}, \text{ for all }r>n+2,
\\
&\Sigma_p  \, \weakstar  \, \Sigma_\infty, \ \ \text{ in }\mM\big(\overline{\Om_T};\R^N \big),
\\
&\sigma_p \, \weakstar \,  \sigma_\infty,  \ \ \ \text{ in }\mM\big(\overline{\Om_T};\R^n \big),
\end{split}
\right.
\eeq
along a subsequence $p_j \to \infty$, where 
\beq
\label{1.29}
\left\{\ \ 
\begin{split}
\Sigma_p\, &:=\, \mathfrak{M}_p \big( \K[u_p] \big)\mL^{n+1}\LL_{\Om_T},
\\
\sigma_p\, &:=\, \mathfrak{M}_p (y_p)\mL^{n+1}\LL_{\Om_T} .
\end{split}
\right.
\eeq
Finally, $\Sigma_\infty$ concentrates on the set whereon $|\K[u_\infty]|$ is maximised over $\overline{\Om_T}$
\beq
\label{1.30}
\Sigma_\infty \Big(\big\{\big|\K[u_\infty]\big| < \big\|\K[u_\infty]\big\|_{L^\infty(\Om_T)} \big\} \Big) =\, 0,
\eeq
and $\sigma_\infty$ asymptotically concentrates on the set whereon $|y_\infty|$ is approximately maximised over $\overline{\Om_T}$, in the sense that for any $\e>0$ small,
\beq
\label{1.31}
\lim_{p\to\infty} \, \sigma_p \Big(\Big\{|y_p| < \|y_\infty\|_{L^\infty(\Om_T)} -\e \Big\} \Big) =\, 0.
\eeq
\end{theorem}

Even though the weak interpretation of the equations \eqref{1.22}-\eqref{1.23} is relatively obvious, this is not the case for \eqref{1.26}-\eqref{1.27} despite having a simpler form. The reason is that the limiting measures $(\Sigma_\infty,\sigma_\infty)$ are not product measures on $\overline{\Om_T}=\overline{\Om}\by [0,T]$ in order to use the Fubini theorem, therefore due to the temporal dependence, \eqref{1.27} cannot be simply interpreted as ``$\div (\si_\infty)=0$". Similar arguments can be made for \eqref{1.26} as well. Since this point is not utilised any further in this paper, we only provide a brief discussion in the next section.

We conclude this introduction with some remarks regarding the organisation of this paper. This introduction is followed by Section \ref{section2}, in which we discuss some preliminaries and also establish some basic estimates which are utilised subsequently to establish our main results. In Section \ref{section3} we prove Theorem \ref{theorem1} by establishing the existence of constrained minimisers for all $p$ including $p=\infty$, as well as the convergence of minimiser of the former problems to those of the latter. In Section \ref{section4} we prove Theorem \ref{theorem2}, deriving the necessary PDE conditions which constrained minimisers in $L^p$ satisfy. Finally, in Section \ref{section5} prove Theorem \ref{theorem3}, deriving the necessary PDE conditions that constrained minimisers in $L^\infty$ satisfy, as well as the additional properties that the measures arising in these PDEs satisfy. A key ingredient here is that we establish appropriate weak* compactness for the Lagrange multipliers arising in the $L^p$ problems in order to pass to the limit as $p\to\infty$.

\section{Preliminaries and the main estimates} \label{section2}

We begin by recording for later use the following modified H\"older inequality for the dotted $\dot L^p$ regularised ``norms" defined in \eqref{1.7}: for any $1\leq q\leq p<\infty$ and $h \in L^p(\Om_T)$, we have the inequality
\[
\| h \|_{\dot L^q (\Om_T)} \, \leq\, \| h \|_{\dot L^p (\Om_T)} \,+ \, \sqrt{q^{-2}-p^{-2}},
\]
which can be very easily confirmed by a direct computation. Next, we continue with a brief discussion regarding the weak interpretation of the equations \eqref{1.26}-\eqref{1.27}. As already noted in the introduction, since $(\Sigma_\infty,\sigma_\infty)$ are not in general neither product measures nor absolutely continuous with respect to the $(n+1)$-Lebesgue measure on $\overline{\Om_T}=\overline{\Om} \by [0,T]$, one needs to use the disintegration ``slicing" theorem for Young measures in order to express them appropriately, as follows. Since $\si_\infty$ is a vector measure in $\mM(\overline{\Om_T};\R^n)$,  by the Radon-Nikodym theorem, we may decompose
\[
\si_\infty \, =\, \frac{\mathrm d \si_\infty}{\mathrm d \|\si_\infty\|} \, \|\si_\infty\| ,
\]
where $ \|\si_\infty\| \in \mM(\overline{\Om_T})$ is the scalar total variation measure and ${\mathrm d \si_\infty}/{\mathrm d \|\si_\infty\|}$ is the vector-valued Radon-Nikodym derivative of $\si_\infty$ with respect to $\|\si_\infty\|$. Fix now any $h \in L^1(\overline{\Om_T},\|\si_\infty\|)$. By the disintegration ``slicing" theorem for Young measure (see se.g.\ \cite[Theorem 3.2, p.\ 179]{FG}), we have the representation formula
\[
\int_{\overline{\Om_T}}h \, \mathrm d \|\si_\infty\| \, =\, \int_{[0,T]}\! \bigg(\int_{\overline{\Om}}h(x, t) \, \mathrm d \|\si_\infty\|_t(x) \bigg) \mathrm d \|\si_\infty\|^o(t)
\] 
where the measure $\|\si_\infty\|^o \in \mM([0,T])$ and the family of measures $(\|\si_\infty\|_t)_{t\in [0,T]} \sub \mM(\overline{\Om})$ are defined as follows:
\[
\|\si_\infty\|^o \,:=\, \|\si_\infty\|\big(\overline{\Om} \by \cdot \, \big), \ \ \ \ \ 
\|\si_\infty\|_t(A) \,:=\, \frac{\mathrm d \|\si_\infty\| \big(A \by \cdot \, \big)}{\mathrm d  \|\si_\infty\| \big(\overline{\Om} \by \cdot \,\big)}(t), \text{ for } A \sub \overline{\Om} \text{ Borel}.
\]
Namely, $\|\si_\infty\|^o$ is one of the marginals of $\si_\infty$ and for $\|\si_\infty\|^o$-a.e.\ $t\in[0,T]$, the measure $\|\si_\infty\|_t$ evaluated at $A$ is defined as the Radon-Nikodym derivative of the measure $\|\si_\infty\|\big(A \by \cdot \, \big)$ with respect to $\|\si_\infty\|\big(\overline{\Om} \by \cdot \, \big)$ at the point $t \in [0,T]$. Then, in view of \eqref{1.27}, by choosing $\mathrm p$ in the form $\mathrm p(x,t)= \pi(x) \tau (t)$, we have
\[
\begin{split}
0 \, & =\, \int_{\overline{\Om_T}} \D \mathrm p \cdot \mathrm d \si_\infty 
\\
&=\, \int_{\overline{\Om_T}} \bigg(\D \mathrm p \cdot \frac{\mathrm d \si_\infty}{\mathrm d \|\si_\infty\|}\bigg) \, \mathrm d \|\si_\infty\|  
\\
&=\, \int_{[0,T]}\! \bigg(\int_{\overline{\Om}}\bigg(\D \mathrm p \cdot \frac{\mathrm d \si_\infty}{\mathrm d \|\si_\infty\|}\bigg)(x, t) \, \mathrm d \|\si_\infty\|_t(x) \bigg) \mathrm d \|\si_\infty\|^o(t)
\\
&=\, \int_{[0,T]}\! \bigg(\int_{\overline{\Om}}\bigg(\D \pi(x) \cdot \frac{\mathrm d \si_\infty}{\mathrm d \|\si_\infty\|}(x, t)\bigg) \, \mathrm d \|\si_\infty\|_t(x) \bigg) \tau(t) \, \mathrm d \|\si_\infty\|^o(t).
\end{split}
\]
The arbitrariness of $\tau$ implies that for $\|\si_\infty\|^o$-a.e.\ $t \in[0,T]$, we have
\[
\int_{\overline{\Om}}\bigg(\D \pi(x) \cdot \frac{\mathrm d \si_\infty}{\mathrm d \|\si_\infty\|}(x, t)\bigg) \, \mathrm d \|\si_\infty\|_t(x) \, =\, 0.
\]
When restricting our attention to those test function for which $\pi|_{\p\Om}\equiv 0$, we obtain the next weak interpretation of \eqref{1.27}:
\[
\div \bigg( \frac{\mathrm d \si_\infty}{\mathrm d \|\si_\infty\|}(\cdot,t) \, \|\si_\infty\|_t \bigg) \, =\, 0, \ \ \ \text{ in }\Om,
\]
for $\|\si_\infty\|^o$-a.e.\ $t \in[0,T]$. Similar considerations apply also to equation \eqref{1.26}, but the arguments are considerably more complicated.

Next we prove a general compact embedding lemma by means of interpolation theory.

\begin{lemma} \label{lemma4} Suppose that $p> n+2$. Then, there exists $\alpha \in (0,1)$ such that the space 
\[
W^{2,1;p}(\Om_T)\,:=\, L^p\big((0,T);W^{2,p}(\Omega)\big) \bigcap W^{1,p}\big((0,T);L^p(\Omega)\big)
\]
is compactly embedded in the space $C^{0,\alpha}\big([0,T];C^{1,\alpha}(\overline{\Omega})\big)$.
\end{lemma}

\BPL \ref{lemma4}. Let us use the abbreviated space notation
\[
\X_1 \, :=\, W^{2,p}(\Omega), \ \  \X_0\, :=\, L^p(\Omega)
\]
and select $\theta$ such that
\[
\frac{p+n}{2p} \, <\, \theta \, <\, \frac{p-1}{p},
\]
which is possible since 
\[
\frac{p-1}{p}-\frac{p+n}{2p} \, =\, \frac{p-(n+2)}{2p} \, >\, 0. 
\]
Since $1-\theta>1/p$, direct application of the interpolation result in \cite[Theorem 5.2]{AmannGlasnik} for the exponents $s_0:=1$, $s_1:=0$ and $p_0\equiv p_1:=p$ yields that space $W^{2,1;p}(\Om_T)$ is compactly embedded in the space $C^{0,\alpha}\big([0,T];\X\big)$, where $0<\alpha<1-\theta-1/p$ and $\X=(\X_0,\X_1)_{\theta,p}$ symbolises the real interpolation between the Banach spaces $\X_0$ and $\X_1$. Now it remains to identify the space $\X$. By using standard results in interpolation theory (see e.g.\ \cite[Theorem 4.3.1.1 and formula (2.4.2/9)]{TriebelInterpolationSpaces} or \cite{Triebel2002} for Lipschitz domains) we get:
\[
\big(L^p(\Omega),W^{2,p}(\Omega)\big)_{\theta,p} \, =\, \mathrm B^{2\theta}_{pp}(\Omega) \, =\, W^{2\theta,p}(\Omega).
\]
Since $2\theta>1+n/p$, by the standard Sobolev embedding theorem for fractional spaces (e.g. \cite[Theorem 8.2]{HitchhikersSobolevSpaces}, we have that $W^{2\theta,p}(\Omega)$ is continuously embedded in the space $C^{1,\alpha}(\overline{\Omega})$, where $0<\alpha\leq  2\theta-1 -n/p$. The conclusion ensues.
\qed
\ms

\begin{remark} \label{remark5} Let us now record for later use the following simple inclusion of space (which is in fact a continuous embedding): 
\[
C^{0,\alpha}\big([0,T];C^{0,\alpha}(\overline{\Omega})\big) \, \sub \, C^{0,\alpha}\big(\overline{\Omega_T}\big).
\]
Indeed, for any $h\in C^{0,\alpha}\big([0,T];C^{0,\alpha}(\overline{\Omega})\big)$, we compute
\[
\begin{split}
\big|h(t_1,x_1)-h(t_2,x_2)\big| \, &\leq \,|h(t_1,x_1)-h(t_2,x_1)|+|h(t_2,x_1)-h(t_2,x_2)|
\\
& \leq \, \|h(t_1, \cdot )-h(t_2, \cdot )\|_{C(\overline{\Omega})} \, +\, \|h(t_2, \cdot )\|_{C^{0,\alpha}(\overline{\Omega})}|x_1-x_2|^{\alpha}
\\
& \leq \, \big(|t_1-t_2|^{\alpha}+|x_1-x_2|^{\alpha}\big)\|h\|_{C^{0,\alpha}([0,T];C^{0,\alpha}(\overline{\Omega}))}
\end{split}
\]
which establishes the claim.
\end{remark}
 
\begin{lemma} Suppose that assumptions \eqref{1.5} and \eqref{1.15} are satisfied.\label{lemma6} We have that 
\[
\mathfrak{X}^\infty(\Om_T) \, \neq \, \emptyset
\]
(and consequently we have $\mathfrak{X}^p(\Om_T) \neq \emptyset$ for all $p>1$). Further, if $(u,\mathrm p, y) \in \mathfrak{X}^p(\Om_T)$ for some $p>1$ which satisfies 
\[
\E_p(u,\mathrm p, y) \, \leq \, M
\]
for some $M>0$, then for any $q\leq p$ there exists $C(q,M)>0$ such that
\[
\|u\|_{W^{2,1;q}_{\mathrm L,\si}(\Om_T)} +\, \|\mathrm p \|_{W^{1,0;q}_{\sharp}(\Om_T)}  +\, \|y \|_{L^q(\Om_T)} \, \leq \, C(q,M).
\]
Further, if $p>n+2$ and $q\in(n+2,p]$, then there exists $\al \in (0,1)$ and a constant $C(M,q)>0$ such that additionally
\[
\|u\|_{C^{0,\al}(\Om_T)} +\, \|\D u\|_{C^{0,\al}(\Om_T)} \, \leq \, C(q,M).
\]
\end{lemma}

We note that the constants above also depend on $n,\p\Om,T,f,u_0,\la$, but as all these are fixed throughout this paper, we suppress denoting the explicit dependence on them.

\BPL \ref{lemma6}. By assumptions \eqref{1.5}(b)-\eqref{1.5}(c), we have that the triplet $(u_0,0,y_0)$, where 
\[
y_0\,:=\, -\nu\Delta u_0 \, +\, (u_0 \cdot \D) u_0 \, - \, f
\]
satisfies that $(u_0,0,y_0) \in \mathfrak{X}^\infty(\Om_T)$, and in fact lies also in the smaller space 
\[
W^{2,1;\infty}_{\mathrm L,\si}(\Om_T) \by \smash{W^{1,0;\infty}_{\sharp}(\Om_T)} \by L^\infty(\Om_T;\R^n).
\]
Next, if $(u,\mathrm p, y) \in \mathfrak{X}^p(\Om_T)$ with $\E_p(u,\mathrm p, y)\leq M$, then we readily have that 
\[
\|y\|_{L^q(\Om_T)} \, \leq \, \|y\|_{\dot L^p(\Om_T)} \, \leq \, \frac{M}{\la}, 
\]
whilst by assumptions \eqref{1.15} and \eqref{1.5}(c) we have that 
\[
\|u\|_{W^{2,1;q}_{\mathrm L,\si}(\Om_T)} +\, \|\mathrm p \|_{W^{1,0;q}_{\sharp}(\Om_T)} \, \leq\, C(q)\Big(1+\frac{M}{\la} \Big).
\]
for some $q$-dependent constant $C(q)$, for any $q\in (n,p]$. Further, suppose $p>n+2$ and $n+2<q\leq p$. Then, the above estimate in particular implies
\[
\|u\|_{L^q(\Om_T)} +\,  \|\nabla u\|_{L^q(\Om_T)}  \, \leq\, C(q,M),
\]
whereat application of the Morrey imbedding theorem yields
\[
\|u\|_{C^{0,\al'}(\Om_T)}  \, \leq\, C(q,M),
\]
for a new constant $C(q,M)$ and some $\al' \in (0,1)$. Next, by Lemma \ref{lemma4}, Remark \ref{remark5} and the established estimate for $q>n+2$, we have
\[
C(q)\|\D u\|_{C^{0,\al''}(\Om_T)} \, \leq\, \|u\|_{W^{2,1;q}(\Om_T)} \, \leq\, C(q,M),
\]
for some $\al'' \in(0,1)$ and some constant $C(q)>0$. By choosing $\al:= \min\{\al',\al''\}$, the conclusion ensues.   \qed
\ms

\section{Minimisers of $L^p$ problems and convergence as $p\to\infty$}  \label{section3}

In this section we establish Theorem \ref{theorem1}, by utilising the results of Section \ref{section2}.

\BPT \ref{theorem1}. Fix $p \in (n+2,\infty)$. By Lemma \ref{lemma4}, $\mathfrak{X}^p(\Om_T)\neq \emptyset$, therefore $0\leq \inf_{\mathfrak{X}^p(\Om_T)} \E_p <\infty$. By Lemma \ref{lemma6}, it follows that $\mathfrak{X}^p(\Om_T)$ is sequentially weakly compact. Note now that $y \mapsto \| y \|^p_{\dot L^p(\Om_T)}$ is trivially convex, and by the identity
\[
\big\| \K\big(\cdot, u,\D u, \p_t u,\mathrm p\big)\big\|^p_{\dot L^p(\Om_T)}\, =\, \int_{\Om_T} \Big(\big|\K\big(\cdot, u,\D u, \p_t u,\mathrm p\big)|^2 + p^{-2}\Big)^{\! \frac{p}{2}}\mathrm d \mL^{n+1},
\]
assumption \eqref{1.5}(f) yields that 
\[
(\p_t u,\mathrm p) \, \mapsto \, \big\| \K\big(\cdot, u,\D u, \p_t u,\mathrm p\big)\big\|^p_{\dot L^p(\Om_T)}
\]
is also convex. Again by Lemma \ref{lemma6} and by standard results in the Calculus of Variations (see e.g.\ \cite{D}), it follows that $\E_p$ is weakly lower semicontinuous in $\mathfrak{X}^p(\Om_T)$ as the convex combination of the $p$-th roots of two weakly lower semicontinuous functionals. Hence, $\E_p$ attains its infimum at some $(u_p,\mathrm p_p , y_p) \in \mathfrak{X}^p(\Om_T)$. 

Consider now the family of minimisers $(u_p,\mathrm p_p , y_p)_{p>n+2}$. For any $(u,\mathrm p, y) \in \mathfrak{X}^\infty(\Om_T)$ and any $q\leq p$, minimality and the H\"older inequality for the dotted $\dot L^p$ functionals yield
\[
\E_p (u_p,\mathrm p_p , y_p) \, \leq\, \E_p (u,\mathrm p , y)\, \leq\, \E_\infty (u,\mathrm p , y) \,+\, p^{-1}.
\]
By choosing $(u,\mathrm p, y)=(u_0,0,y_0)$, by Lemma \ref{lemma6} and a standard diagonal argument, we have that the family of minimisers is weakly precompact in $\mathcal{W}^q(\Om_T)$ for all $q \in (n+2,\infty)$. Further, by Lemma \ref{lemma4} and Remark \ref{remark5}, $W^{2,1;q}_{\mathrm L,\si}(\Om_T;\R^n)$ is compactly embedded in $\smash{C^{0,\al}\big([0,T];C^{1,\al}(\overline{\Om};\R^n)\big)}$. Hence, for any sequence of indices $p_j \to \infty$, there exists $(u_\infty,\mathrm p_\infty , y_\infty) \in \cap_{q\in(n+2,\infty)}\mathcal{W}^q(\Om_T)$ and a subsequence denoted again as $(p_j)_1^\infty$ such that \eqref{1.17} holds true. Additionally, due to these modes of convergence, it follows that $(u_\infty,\mathrm p_\infty , y_\infty)$ solves \eqref{1.1}, therefore in fact $(u_\infty,\mathrm p_\infty , y_\infty) \in \mathfrak{X}^\infty(\Om_T)$. Again now by minimality and the H\"older inequality for the dotted $\dot L^p$ functionals, for any $(u,\mathrm p, y) \in \mathfrak{X}^\infty(\Om_T)$ we have
\[
\E_q (u_p,\mathrm p_p , y_p) -\sqrt{q^{-2}-p^{-2}}\, \leq\, \E_p (u_p,\mathrm p_p , y_p) \, \leq\, \E_\infty (u,\mathrm p , y) \,+\, p^{-1}.
\]
Since as we have already shown, $\E_q$ is weakly lower semicontinuous in $\mathfrak{X}^q(\Om_T)$, by letting $p \to \infty$ along the subsequence in the above inequality yields
\[
\begin{split}
\E_q (u_\infty,\mathrm p_\infty , y_\infty) -\sqrt{q^{-2}}\, 
&\leq\, \liminf_{p_j\to \infty}\E_p (u_p,\mathrm p_p , y_p) 
\\
& \leq\, \limsup_{p_j\to \infty}\E_p (u_p,\mathrm p_p , y_p) 
\\
& \leq\, \E_\infty (u,\mathrm p , y) .
\end{split}
\]
By further letting $q\to \infty$, we obtain
\[
\begin{split}
\E_\infty (u_\infty,\mathrm p_\infty , y_\infty) \, &\leq\, \liminf_{p_j\to \infty}\E_p (u_p,\mathrm p_p , y_p) 
\\
& \leq\, \limsup_{p_j\to \infty}\E_p (u_p,\mathrm p_p , y_p) 
\\
& \leq\, \E_\infty (u,\mathrm p , y) ,
\end{split}
\]
for any $(u,\mathrm p, y) \in \mathfrak{X}^\infty(\Om_T)$. The above inequality establishes on the one hand that $(u_\infty,\mathrm p_\infty , y_\infty)$ minimises $\E_\infty$ over $\mathfrak{X}^\infty(\Om_T)$, and on the other hand by choosing $(u ,\mathrm p , y ):=(u_\infty,\mathrm p_\infty , y_\infty)$ that \eqref{1.18} holds true. Hence, Theorem \ref{theorem1} has been established.
\qed
\ms

\section{The equations for $L^p$ PDE-constrained minimisers} \label{section4}

In this section we establish the proof of Theorem \ref{theorem2}. We begin with some preparation. Firstly, it will be convenient to consider the functional $\E_p$ of \eqref{1.6} as being defined in the wider Banach space $\mathcal W^p(\Om_T)$ defined in  \eqref{1.9}:
\[
\E_p \ : \ \ \mathcal W^p(\Om_T) \larrow \R.
\]
Next, we introduce a mapping on $\mathcal W^p(\Om_T)$ which incorporates the PDE constraint \eqref{1.1} appearing in \eqref{1.8} as follows. We define
\[
\G=\left[
\begin{array}{c}\G_1 
\\
\G_2
\end{array} 
\right]
 \ : \ \ \mathcal W^p(\Om_T) \larrow L^p(\Om_T;\R^n)\times W^{2-\frac{2}{p},p}_{0,\si}(\Om;\R^n)
\]
by setting
\[
\left\{ \ \
\begin{split}
\G_1(u,\mathrm p, y) \, & := \, \p_t u \,-\nu\Delta u\,+ (u\cdot \D)u\, +\D \mathrm p \, -(y+f),
\\
\G_2(u,\mathrm p, y) \, &: =\, u(\cdot,0)-u_0.
\end{split}
\right.
\]
Then, we may express \eqref{1.8} as
\[
\mX^p(\Om_T) \, =\, \mathcal W^p(\Om_T) \cap \{\G=0\}.
\]
We are now ready to prove our second main result.

\BPT \ref{theorem2}. By assumption \eqref{1.5}, for any $p\in (n+2,\infty)$ the functional $\E_p : \mathcal W^p(\Om_T) \larrow \R$ is Frech\'et differentiable and its derivative
\[
\begin{split}
\mathrm d  \E_p  \ & :  \ \ \mathcal W^p (\Om_T) \larrow \big( \mathcal W^p(\Om_T)\big)^*, 
\\
 (\mathrm d \E_p )_{ (\bar u , \bar{\mathrm p}, \bar y)} (u , \mathrm p, y) \, & =\, \frac{\mathrm d}{\mathrm d \e}\Big|_{\e=0} \E_p\Big( \bar u +\e u, \bar{\mathrm p}+\e \mathrm p, \bar y +\e y\Big)
\end{split}
\]
can be easily computed and is given by the formula
\[
\begin{split}
(\mathrm d \E_p )_{ (\bar u , \bar{\mathrm p}, \bar y)} (u , \mathrm p, y) \, =\, p(1-\lambda)\,  & {\av_{\Om_T}}  \! \Big(\K_\eta [\bar u,\bar {\mathrm p}]\cdot u \, +\, \K_{(A,a)}[\bar u,\bar {\mathrm p}] : \nabla u +\,  \K_r[\bar u,\bar{\mathrm p}] \, \mathrm p\Big) \cdot
\\
& \cdot \mathfrak{M}_p \big( \K[\bar u,\bar{\mathrm p}] \big) \, \d\mL^{n+1}
\, +\, p\la \,{\av_{\Om_T}} \mathfrak{M}_p(\bar y)\cdot y \, \d \mL^{n+1},
\end{split}
\]
where the operator $\mathfrak{M}_p : L^p(\Om_T;\R^M) \larrow L^{p'}(\Om_T;\R^M)$ (for $M\in\{N,n\}$) is given by \eqref{1.19} and we have used the notation introduced in \eqref{K[u,p]}. Next, we note that the mapping $\G$ which incorporates the PDE constraint is also Fr\'echet differentiable and it can be easily confirmed that its derivative
\[
\begin{split}
\mathrm d \G  \ &  : \ \ \mathcal W^p(\Om_T) \larrow \mathcal B\bigg(\mathcal W^p (\Om_T), L^p(\Om_T;\R^n)\times W^{2-\frac{2}{p},p}_{0,\si}(\Om;\R^n)\bigg),
\\
& (\mathrm d \G )_{ (\bar u , \bar{\mathrm p}, \bar y)} (u , \mathrm p, y) \,  =\, \frac{\mathrm d}{\mathrm d \e}\Big|_{\e=0} \G \Big( \bar u +\e u, \bar{\mathrm p}+\e \mathrm p, \bar y +\e y\Big)
\end{split}
\]
is given by the formula
\[
\begin{split}
 \big(\mathrm d   \G\big)_{ (\bar u , \bar{\mathrm p}, \bar y)} (u , \mathrm p, y)\, = \, 
\left[
\begin{array}{c}
\p_t u -\nu \Delta u + (u\cdot \D)\bar{u} + (\bar{u}\cdot \D)u + \D  \mathrm p - y
 \ms
 \\
u(\cdot,0)
 \end{array}
\right] .
\end{split}
\]
We now claim that the differential
\[
(\mathrm d \G )_{ (\bar u , \bar{\mathrm p}, \bar y)} \ : \ \ \mathcal W^p (\Om_T) \larrow L^p(\Om_T;\R^n)\times W^{2-\frac{2}{p},p}_{0,\si}(\Om;\R^n)
\]
is a surjective map, for any $(\bar u , \bar{\mathrm p}, \bar y) \in \mathcal W^p(\Om)$. This is equivalent to the statement that for any $p>n+2$, the linearised Navier-Stokes problem
\[
\left\{ \ \
\begin{array}{rr}
\p_t u -\nu \Delta u + (u\cdot \D)\bar{u} + (\bar{u}\cdot \D)u + \D  \mathrm p \, =\, F, & \text{ in }\Om_T,
\\
\div \, u \, =\, 0, & \text{ in }\Om_T,
\\
u(\cdot,0) \, =\, v, & \text{ on }\Om,
\\
u \, =\, 0, & \text{ on }\p \Om \by (0,T),
\end{array}
\right.
\]
has a solution $(u,\mathrm p) \in W^{2,1;p}_{\mathrm L,\si}(\Om_T;\R^n) \by W^{1,0;p}_\sharp(\Om_T)$, for any $\bar u \in W^{2,1;p}_{\mathrm L,\si}(\Om_T;\R^n)$ and any data
\[
(F,v) \, \in \, L^p(\Om_T;\R^n)\times W^{2-\frac{2}{p},p}_{0,\si}(\Om;\R^n).
\]
This is indeed the case, and it is a consequence of a classical result of Solonnikov \cite[Th.\ 4.2]{S} for $n=3$ and of Giga-Sohr \cite[Th.\ 2.8]{GS} for $n>3$, as a perturbation of the Stokes problem. As a consequence, the assumptions of the generalised Kuhn-Tucker theorem hold true (see e.g.\ Zeidler \cite[Cor.\ 48.10 \& Th.\ 48B]{Z}). Hence, there exists a Lagrange multiplier
\[
\begin{split}
\Lambda_p \ \in \ & \Big(L^{p} (\Om_T;\R^{n}) \by W^{2-\frac{2}{p},p}_{0,\si} (\Om_T;\R^n)\Big)^*
\end{split}
\]
such that
\[
\begin{split}
 \big(\mathrm d   \E_p\big)_{ (u _p, \mathrm p_p,y_p)}   (u,\mathrm p,y) \, =\, \Big\langle \big(\mathrm d   \G\big)_{ (u _p,\mathrm p_p,y_p)}(u,\mathrm p,y),\,  \Lambda_p  \Big\rangle ,
\\
\end{split}
\]
for any $(u , {\mathrm p}, y) \in \mathcal W^p(\Om)$. By standard duality arguments, the Riesz representation theorem and by taking into account the form  of the differentials $\mathrm d \E_p$ and $\mathrm d \G$, we may identify $\La_p$ with a pair of Lagrange multipliers
\[
 (\phi_p, \Psi_p  ) \, \in \, L^{p'} (\Om_T;\R^n) \by \Big( W^{2-\frac{2}{p},p}_{0,\si} (\Om_T;\R^n)\Big)^*
\]
such that, the constrained minimiser $\big(u _p, p_p, y_p\big) \in \mX^p(\Om_T)$ satisfies the equation
 \[
\begin{split}
& (1-\lambda)  \int_{\Om_T}  \! \Big(\K_\eta [u_p,{\mathrm p}_p]\cdot u \, +\, \K_{(A,a)}[u_p,{\mathrm p}_p] : \nabla u +\,  \K_r[u_p,{\mathrm p}_p] \, \mathrm p\Big)\cdot
\\
&\hspace{50pt}\cdot \mathfrak{M}_p \big( \K[u_p,{\mathrm p}_p] \big) \, \d\mL^{n+1}  \ +\ \la \int_{\Om_T} \mathfrak{M}_p(y_p)\cdot y \, \d \mL^{n+1}
\\
& = \int_{\Om_T} \Big(\p_t u -\nu \Delta u + (u\cdot \D)u_p + (u_p \cdot \D)u + \D  \mathrm p - y \Big) \cdot \phi_p \, \d \mL^{n+1} \,+\, \langle \Psi_p, u(\cdot,0)\rangle ,
\end{split}
\]
for any $(u,\mathrm p,y) \in \mathcal W^p(\Om_T)$. We note that here we have tacitly rescaled $(\phi_p, \Psi_p  )$ by multiplying them with the factor $p(\mL^{n+1}(\Om_T))^{-1}$, in order to remove the averages arising from $\mathrm \E_p$ on the left hand side and to be able to obtain non-trivial limits as $p\to \infty$ of the multipliers themselves later on. By using linear independence, the above equation actually decouples to the triplet of relations
\[
\left\{ \ \
\begin{split}
& (1-\lambda)  \int_{\Om_T}  \! \Big(\K_\eta [u_p,{\mathrm p}_p]\cdot u \, +\, \K_{(A,a)}[u_p,{\mathrm p}_p] : \nabla u \Big) \cdot \mathfrak{M}_p \big( \K[u_p,{\mathrm p}_p] \big) \, \d\mL^{n+1}  \ 
\\
& = \int_{\Om_T} \Big(\p_t u -\nu \Delta u + (u\cdot \D)u_p + (u_p \cdot \D)u \Big) \cdot \phi_p \, \d \mL^{n+1} \,+\, \langle \Psi_p, u(\cdot,0)\rangle ,
\end{split}
\right.
\]
\[
\begin{split}
 (1-\lambda)  \int_{\Om_T}  \! \big( \K_r[u_p,{\mathrm p}_p] \, \mathrm p\big) \cdot \mathfrak{M}_p \big( \K[u_p,{\mathrm p}_p] \big) \, \d\mL^{n+1}  
\, = \, \int_{\Om_T}   \D  \mathrm p   \cdot \phi_p \, \d \mL^{n+1},
\end{split}
\]
\[
\begin{split}
 \la \int_{\Om_T} \mathfrak{M}_p(y_p)\cdot y \, \d \mL^{n+1}
\, = \, -\int_{\Om_T} y  \cdot \phi_p \, \d \mL^{n+1} .
\end{split}
\]
The arbitrariness of $y \in L^p(\Om_T ;\R^n)$ in the third equation readily yields that the multiplier $\phi_p$ equals
\[
\phi_p \,=\, -\la \,\mathfrak{M}_p(y_p).
\] 
By substituting this into the first two equations, we see that the theorem has been established. \qed
\ms

\section{The equations for $L^\infty$ PDE-constrained minimisers} \label{section5}
 
In this section we establish our final main result. In this case we need to assume that $\K$ does not depend on $(\p_t u, \mathrm p)$ and we will invoke the symbolisations \eqref{K[u]} for $\K$ and its derivatives $\K_\eta,\K_A$, all of which are additionally assumed to be continuous.

\BPT \ref{theorem3}. By Theorem \ref{theorem2} it follows that for any $p \in (n+2,\infty)$, the minimising triplet $(u_p,\mathrm p_p, y_p)\in\mathfrak{X}^p(\Om_T)$ satisfies
\[
\left\{ \ 
\begin{split}
 &(1-\lambda)\int_{\Om_T} \Big(\K_\eta [u_p]\cdot u \, +\, \K_{A}[u_p] : \D u \Big)\cdot \mathfrak{M}_p \big( \K[u_p] \big) \, \d\mL^{n+1}
 \\
 &= 
 -\la  \int _{\Om_T}\Big(\p_t u - \nu\Delta u + (u\cdot \D)u_p+( u_p\cdot \D)u \Big)\cdot \mathfrak{M}_p(y_p) \, \d {\mathcal{L}}^{n+1}\, + \, \big\langle \Psi_p , u(\cdot,0) \big\rangle,
\end{split}
\right.
\]
and also
\[
\begin{split}
& \int_{\Om_T}\D \mathrm p \cdot \mathfrak{M}_p(y_p) \, \d \mL^{n+1} \, =\, 0,
\\
\end{split}
\]
for all test mappings $(u, \mathrm p)\in W^{2,1;p}_{\mathrm L,\si}(\Om_T;\R^n) \by W^{1,0;p}_\sharp (\Om_T)$. The first goal is to pass to the limit as $p\to \infty$ in these equations in order to obtain \eqref{1.26}-\eqref{1.27}. Since by \eqref{1.19} we readily have that $\mathfrak{M}_p(y_p)$ and $\mathfrak{M}_p\big(\K[u_p]\big)$ are valued in the unit balls of $\smash{L^{p'}}(\Om_T;\R^n)$ and of $\smash{L^{p'}}(\Om_T;\R^N)$ respectively, by defining $\Sigma_p$ and $\sigma_p$ as in \eqref{1.29}, the existence of limiting measures $\Si_\infty,\si_\infty$ is guaranteed along perhaps a further subsequence such that 
\[
\text{$\Sigma_p \weakstar  \Sigma_\infty$ in $\mM\big(\overline{\Om_T};\R^N \big)$ \ and \ $\sigma_p \weakstar  \sigma_\infty$ in $\mM\big(\overline{\Om_T};\R^n \big)$,} 
\]
as $p_j\to \infty$. Further, by Lemma \ref{lemma6} we have that $u_p \larrow u_\infty$ and $\D u_p \larrow \D u_\infty$, both uniformly on $\overline{\Om_T}$ as $p_j\to \infty$. Also, by the continuity assumption on $\K,\K_\eta,\K_A$ on $\overline{\Om_T} \by \R^n \by \R^{n\by n}$ and again Lemma \ref{lemma6}, it follows that 
\[
\text{$\K[u_p] \larrow \K[u_\infty]$, \ $\K_\eta[u_p] \larrow \K_\eta[u_\infty]$ \ and \ $\K_A[u_p] \larrow \K_A[u_\infty]$,} 
\]
all uniformly on $\overline{\Om_T}$ as $p_j\to \infty$. Putting all this together, we see that the remaining main point is to obtain a uniform estimate on the family of Lagrange multipliers $(\Psi_p)_{p>n+2}$ in order to deduce that 
\[
\text{$\Psi_p  \weakstar   \Psi_\infty$ \ in  $\smash{\big(W^{2-2/r,r}_{0,\si} (\Om;\R^{n})\big)^*}$, \ for all $r>n+2$,} 
\]
which would allow to pass to the limit as $p_j\to \infty$. Once this has been achieved, passing to the limit in the equations follows by standard duality pairing arguments, which are made possible by restricting the class of test functions $(u, \mathrm p)$ to those which are continuous together with those  derivatives appearing in the relations.

In order to derive the desired estimate on $(\Psi_p)_{p >n+2}$, we argue as follows. Consider \eqref{1.22} for $\K_a \equiv 0$ (the first equation appearing in this proof) and let us fix the initial value on $\Om\by \{0\}$
\[
u(\cdot,0) \, \equiv \, \hat u \, \in W^{2,\infty}_{0,\si}(\Om;\R^n)
\]
of the arbitrary test function $u$, but we will select $u$ on $\Om_T$ such that the term in the bracket in the integral on the right-hand-side becomes a gradient. Then, this term will vanish identically as a consequence of \eqref{1.23} when $\K_r \equiv 0$ (the second equation appearing in this proof). Indeed, let $p>n+2$ and let also $(\tilde u, \tilde{\mathrm p})$ be the (unique) solution to
\[
\left\{ \ \
\begin{array}{rr}
\p_t \tilde u -\nu \Delta \tilde u + (\tilde u\cdot \D)u_p + (u_p \cdot \D)\tilde u + \D  \tilde{\mathrm p} \, =\, 0, & \text{ in }\Om_T,
\\
\div \, \tilde u \, =\, 0, & \text{ in }\Om_T,
\\
\tilde u(\cdot,0) \, =\, \hat u, & \text{ on }\Om,
\\
\tilde u \, =\, 0, & \text{ on }\p \Om \by (0,T),
\end{array}
\right.
\]
The solvability of the above problem is a consequence of the classical result of Solonnikov \cite[Th.\ 4.2]{S} for $n=3$ and of Giga-Sohr \cite[Th.\ 2.8]{GS} for $n>3$, as a perturbation of the Stokes problem: by choosing $q>n+2$ in Solonnikov's assumption (4.14), a solution as claimed does exist. Further, since $\hat u$ is in $W^{2,\infty}_{0,\si}(\Om;\R^n)$, by \cite[Cor.\ 2, p.\ 489]{S} we have the uniform estimate
\[
\| \tilde u\|_{W^{2,1;r}_{\mathrm L,\si}(\Om_T)} + \, \| \tilde{\mathrm p} \|_{W^{1,0;r}_\sharp (\Om_T)} \, \leq\, C(r) \| \hat u \|_{W^{2-\frac{2}{r},r}_{0,\si}(\Om)},
\] 
for any $r\in (1,\infty)$. By Lemmas \ref{lemma4} and \ref{lemma6} and Remark \ref{remark5}, if we restrict our attention to $r\in (n+2,\infty)$, we additionally have the bound
\[
\| \tilde u\|_{L^\infty(\Om_T)} +\, \| \D \tilde u\|_{L^\infty(\Om_T)}  \leq\, C(r) \| \hat u \|_{W^{2-\frac{2}{r},r}_{0,\si}(\Om)},
\] 
for some new constant $C(r)$ (which is unbounded as $r\searrow n+2$). By setting
\[
\left\{ \ \ 
\begin{split}
K_\infty \,:=& \, \sup \Big\{|\K_\eta|+|\K_A| : \ \Om_T \by \mB^n_{R_\infty}(0) \by \mB^{n\by n}_{R_\infty}(0) \Big\}, 
\\
R_\infty\,:=&\, \sup_{p>n+2}\Big(\| u_p \|_{L^\infty(\Om_T)} +\,\| \D u_p \|_{L^\infty(\Om_T)} \Big),
\end{split}
\right.
\]
where $\mB^n_{R_\infty}(0)$ and $\mB^{n\by n}_{R_\infty}(0)$ denote the balls of radius $R_\infty$ centred at the origin of $\R^n$ and of $\R^{n\by n}$ respectively, we estimate by using \eqref{1.22}-\eqref{1.23} (for $\K_a\equiv 0$, $\K_r\equiv 0$) and that by \eqref{1.19} we have $\smash{\big\| \mathfrak{M}_p \big( \K[u_p] \big)\big \|_{L^1(\Om_T)}\leq 1}$ (for the normalised $L^1$ norm):
\[
\begin{split}
\big| \big\langle \Psi_p , \hat u \big\rangle \big| \,  \leq &\ \la \,  \bigg|\! \int _{\Om_T} \D \tilde{\mathrm p}\cdot \mathfrak{M}_p(y_p) \, \d {\mathcal{L}}^{n+1} \bigg|
\\
&  + (1-\lambda)\bigg| \! \int_{\Om_T} \! \Big(\K_\eta [u_p]\cdot \tilde u \, +\, \K_{A}[u_p] : \D \tilde u \Big)\cdot \mathfrak{M}_p \big( \K[u_p] \big) \, \d\mL^{n+1}\bigg| 
\\
= &\ (1-\lambda)T \mL^n(\Om)\bigg| \,  {\av_{\Om_T}} \! \Big(\K_\eta [u_p]\cdot \tilde u \, +\, \K_{A}[u_p] : \D \tilde u \Big)\cdot \mathfrak{M}_p \big( \K[u_p] \big) \, \d\mL^{n+1}\bigg| 
\\
\leq &\ (1-\la)T \mL^n(\Om) K_\infty   \Big(\| \tilde u \|_{L^\infty(\Om_T)} +\,\| \D \tilde u \|_{L^\infty(\Om_T)} \Big)
\\
\leq &\ (1-\la)T \mL^n(\Om) K_\infty C(r) \| \hat u \|_{W^{2-\frac{2}{r},r}_{0,\si}(\Om)},
\end{split}
\]
for any $r$ fixed. Therefore, for any $r\in (n+2,\infty)$ and any $\hat u \in W^{2,\infty}_{0,\si} (\Om;\R^n)$, we have the estimate
\[
\begin{split}
\big| \big\langle \Psi_p , \hat u \big\rangle \big| \, \leq \big( T \mL^n(\Om) K_\infty  \big)C(r) \| \hat u \|_{W^{2-\frac{2}{r},r}_{0,\si}(\Om)}.
\end{split}
\]
Since $W^{2,\infty}_{0,\si} (\Om;\R^n)$ is dense in $W^{2-2/r,r}_{0,\si} (\Om;\R^n)$, by the Hahn-Banach theorem, the above estimate implies that for any fixed $p>n+2$, the bounded linear functional
\[
\Psi_p \ :\ \  W^{2-\frac{2}{p},p}_{0,\si} (\Om;\R^n) \larrow \R
\]
can be (uniquely) extended to a functional $\Psi_p : W^{2-2/r,r}_{0,\si} (\Om;\R^n) \larrow \R$ for all $r \in (n+2, p]$, whose extension we denote again by $\Psi_p$. Therefore, $\Psi_p$ can be extended to a unique continuous linear functional
\[
\Psi_p \ :\ \  \bigcup_{r>n+2} W^{2-\frac{2}{r},r}_{0,\si} (\Om;\R^n) \larrow \R
\]
on the above Fr\'echet space, whose topology can be defined in the standard locally convex sense by the family of seminorms
\[
\Big\{ \|\cdot\, \|_{W^{2-{2}/{r},r}(\Om)} :\ r>n+2\Big\}. 
\]
Additionally, the uniformity of the estimate with respect to $p$ implies that
\[
\text{$(\Psi_p)_{p>n+2}$ is bounded in $\bigg( \bigcup_{r>n+2}W^{2-\frac{2}{r},r}_{0,\si} (\Om;\R^n) \bigg)^*$}
\]
(in the locally convex sense). Hence, as it can be seen by a customary diagonal argument in the scale of Banach spaces $\smash{\big\{W^{2-{2}/{r},r}_{0,\si} (\Om_T;\R^n) : r>n+2 \big\}}$ comprising the Fr\'echet space, there exists a continuous linear functional
\[
\Psi_\infty  \ :\ \  \bigcup_{r>n+2} W^{2-\frac{2}{r},r}_{0,\si} (\Om;\R^n) \larrow \R
\]
and a further subsequence as $p\to\infty$ such that along which we have $\Psi_p \weakstar \Psi_\infty$ in the locally convex sense. Additionally, since
\[
\Psi_\infty  \in   \bigcap_{r>n+2} \Big(W^{2-\frac{2}{r},r}_{0,\si} (\Om;\R^n)\Big)^* 
\]
the convergence $\Psi_p \weakstar \Psi_\infty$ is equivalent to weak* convergence in the  Banach space $\smash{W^{2-{2}/{r},r}_{0,\si} (\Om;\R^n)}$ for any fixed $r>n+2$. In conclusion, we see that \eqref{1.26}-\eqref{1.27} have now been established.

Now we complete the proof of Theorem \ref{theorem3} by establishing \eqref{1.30}-\eqref{1.31}. Since $\K[u_p]\larrow \K[u_\infty]$ in $C\big(\overline{\Om_T};\R^N\big)$, by applying \cite[Prop.\ 10]{K2}, we immediately obtain that $\Sigma_\infty$ concentrates on the set whereon $|\K[u_\infty]|$ is maximised over $\overline{\Om_T}$:
\[
\Sigma_\infty \bigg(\Big\{\big|\K[u_\infty]\big| < \max_{\overline{\Om_T}}\big|\K[u_\infty]\big| \Big\} \bigg) =\, 0.
\]
This proves \eqref{1.30}. For \eqref{1.31}, we argue as follows. We first note that
\[
\| y_p \|_{\dot L^p(\Om_T)} \larrow \, \| y_\infty \|_{L^\infty(\Om_T)}  
\]
as $p\to \infty$, along a subsequence. In view of \eqref{1.6} and \eqref{1.13}, this is a consequence of \eqref{1.18} and the fact that $\K[u_p]\larrow \K[u_\infty]$ uniformly on $\overline{\Om_T}$, which implies
\[
\big\|\K[u_p] \big\|_{\dot L^p(\Om_T)} \larrow \, \big\| \K[u_\infty] \big\|_{L^\infty(\Om_T)}.  
\]
As a consequence of the convergence of $\| y_p \|_{\dot L^p(\Om_T)}$ to $\| y_\infty \|_{L^\infty(\Om_T)}$, for any $\e>0$ we may choose $p$ large so that 
\[
\| y_p \|_{\dot L^p(\Om_T)} \geq\, \| y_\infty \|_{L^\infty(\Om_T)} -\frac{\e}{2}.
\]
Let us define now the following subset of $\Om_T$, which without loss of generality we may assume it has positive $\mL^{n+1}$-measure:
\[
A_{p,\e} \, :=\, \Big\{ |y_p| \leq \| y_\infty \|_{L^\infty(\Om_T)} -\e\Big\}.
\]
In particular, if $\mL^{n+1}( A_{p,\e})>0$, then necessarily $\| y_\infty \|_{L^\infty(\Om_T)}>0$. For any Borel set $B \sub \Om_T$ such that $\mL^{n+1}(\Om_T \cap B)>0$, we estimate by using \eqref{1.29}, \eqref{1.19}, \eqref{1.7} and the above:
\[
\begin{split}
\si_p( A_{p,\e} \cap B) \, & \leq \, \frac{\mL^{n+1}( A_{p,\e} \cap B)}{\| y_p \|^{p-1}_{\dot L^p(\Om_T)}} \, {\av_{ A_{p,\e} \cap B}} \big(|y_p|_{(p)}\big)^{p-1}\, \mathrm d \mL^{n+1}
\\
& \leq \, \frac{\mL^{n+1}( A_{p,\e} \cap B)}{\| y_p \|^{p-1}_{\dot L^p(\Om_T)}} \, {\av_{ A_{p,\e} \cap B}} \Big( \| y_\infty \|_{L^\infty(\Om_T)} -\e\Big)^{\! p-1}\, \mathrm d \mL^{n+1}
\\
 &\leq \, \frac{\mL^{n+1}( A_{p,\e} \cap B)}{\| y_p \|^{p-1}_{\dot L^p(\Om_T)}}  \Big( \| y_\infty \|_{L^\infty(\Om_T)} -\e\Big)^{\! p-1}
\\
&\leq \, \mL^{n+1}( A_{p,\e} \cap B) \bigg(\frac{\| y_\infty \|_{L^\infty(\Om_T)} -\e}{\| y_\infty \|_{L^\infty(\Om_T)} -\frac{\e}{2}} \bigg)^{\! p-1}.
\end{split}
\]
As a result, for any $\e>0$ small, any $p$ large enough and any Borel set $B \sub \Om_T$ with $\mL^{n+1}(\Om_T \cap B)>0$, we have obtained the density estimate
\[
\begin{split}
\frac{\si_p( A_{p,\e} \cap B)}{\mL^{n+1}( A_{p,\e} \cap B)} \, \leq \, \bigg( 1 - \frac{\e}{2\| y_\infty \|_{L^\infty(\Om_T)} -\e} \bigg)^{\! p-1}.
\end{split}
\]
The above estimate in particular implies that $\sigma_p ( A_{p,\e}) \larrow 0$ as $p\to\infty$ for any $\e>0$ fixed, therefore establishing \eqref{1.31}. The proof of Theorem \ref{theorem3} is now complete.  \qed
\ms

\begin{remark} It is perhaps worth noting (in relation to the preceding arguments in the proof of \eqref{1.31}) that the modes of convergence
\[
\| y_p \|_{L^p(\Om_T)} \larrow \, \| y_\infty \|_{L^\infty(\Om_T)} \ \ \text{ and }\ \ y_p \weakstar y_\infty \ \text{ in }L^\infty(\Om_T;\R^n)
\]
as $p\to \infty$, in general by themselves do not suffice to obtain $y_p \larrow y_\infty$ in any strong sense, hence precluding the derivation of a stronger property than \eqref{1.31}, along the lines of \eqref{1.30}. A simple counter-example, even in one dimension, is the following: let $p\in2\N$ and set
\[
y_p \,:=\, \sum_{j=0}^{(p-2)/2} \bigg[ \chi_{\big(\frac{2j}{p},\frac{2j+1}{p}\big)} - \chi_{\big(\frac{2j+1}{p},\frac{2j+2}{p}\big)}\bigg] \, +\, \chi_{(1,2)},
\]
and also $y_\infty:=\chi_{(1,2)}$. Then, we have $|y_p|=1$ $\mL^1$-a.e.\ on $(0,2)$ for all $p$, hence we deduce that $\| y_p \|_{L^p(0,2)} \larrow \, \| y_\infty \|_{L^\infty(0,2)}$, whilst we also have $\ y_p \weakstar y_\infty$ in $L^\infty(0,2)$ as $p\to \infty$, but $y_p \,\, \not\!\!\larrow y_\infty$ neither a.e., nor in $L^1$ or in measure.
\end{remark}
\ms

\subsection*{Acknowledgement} N.K.\ is indebted to Jochen Br\"ocker for his expert insights involving scientific discussions about variational data assimilation in continuous time in relation to meteorology.

\ms

\bibliographystyle{amsplain}

\begin{thebibliography}{30}


\bibitem{AmannGlasnik}
H. Amann,
\newblock Compact embeddings of vector-valued Sobolev and Besov spaces,
\newblock volume 35(55), pages 161--177, 2000 
\newblock (Dedicated to the memory of Branko Najman).

\bibitem{Amann2} H. Amann, \emph{On the Strong Solvability of the Navier-Stokes Equations}, Journal of Mathematical Fluid Mechanics 2, 16-98 (2000). 

\bibitem{AP} N. Ansini, F. Prinari, \emph{On the lower semicontinuity of supremal functional under differential constraints}, ESAIM - Control, Opt. and Calc. Var. 21(4), 1053-1075 (2015).

\bibitem{A1} G. Aronsson, \emph{Minimization problems for the functional $sup_x \mF(x,
f(x), f'(x))$}, Arkiv f\"ur Mat. 6 (1965), 33 - 53.

\bibitem{A2} G. Aronsson, \emph{Minimization problems for the functional $sup_x \mF(x,
f(x), f'(x))$ II}, Arkiv f\"ur Mat. 6 (1966), 409 - 431.
 
\bibitem{AB} G. Aronsson, E.N. Barron, \emph{$L^\infty$ Variational Problems with Running Costs and Constraints}, Appl. Math. Optimization 65, 53-90 (2012).

\bibitem{AK1} B. Ayanbayev, N. Katzourakis, \emph{Vectorial variational principles in $L^\infty$ and their characterisation through PDE systems}, Applied Mathematics \& Optimization, 1-16 (2019).

\bibitem{AK2} B. Ayanbayev, N. Katzourakis, \emph{A Pointwise Characterisation of the PDE system of vectorial Calculus of variations in $L^\infty$}, Proceedings of the Royal Society of Edinburgh: Section A Mathematics, 1-17. doi:10.1017/prm.2018.89.

\bibitem{BBJ} E.N. Barron, M. Bocea, R. Jensen, \emph{Viscosity solutions of stationary Hamilton-Jacobi equations and minimizers of $L^\infty$ functionals}, Proc. Amer. Math. 145(12), 5257-5265 (2017).

\bibitem{BJ} E.N. Barron, R. Jensen, \emph{Minimizing the $L^\infty$ norm of the gradient with an energy constraint}, Comm. Partial Differential Equations 30, 10-12, 1741-1772 (2005).

\bibitem{BJW1} E. N. Barron, R. Jensen, C. Wang, \emph{The Euler equation and absolute minimizers of $L^{\infty}$ functionals}, Arch. Rational Mech. Analysis 157 (2001), 255-283.

\bibitem{BJW2} E. N. Barron, R. Jensen, C. Wang, \emph{Lower Semicontinuity of $L^{\infty}$ Functionals} Ann. I. H. Poincar\'e AN 18, 4 (2001)
495-517.

\bibitem{BOT} H. Bessail, E. Olson, E.S. Titi, \emph{Continuous data assimilation with stochastically noisy data}, Nonlinearity 28 (2015) 729-753.

\bibitem{BN} M. Bocea, V. Nesi, \emph{$\Ga$-convergence of power-law functionals, variational principles in $L^\infty$, and applications}, SIAM J. Math. Anal., 39 (2008), 1550-1576.

\bibitem{BP} M. Bocea, C. Popovici, \emph{Variational principles in $L^\infty$ with applications to antiplane shear and plane stress plasticity}, Journal of Convex Analysis Vol. 18 No. 2, (2011) 403-416.

\bibitem{B1} J. Br\"ocker, \emph{What is the correct cost functional for variational data assimilation?} Climate Dynamics, 52 (1-2), 389-399 (2019). 

\bibitem{B2} J. Br\"oecker, \emph{On variational data assimilation in continuous time}, Quarterly Journal of the Royal Meteorological Society 136:652, Part A, 1906-1919 (2010).

\bibitem{B3} J. Br\"oecker, \emph{Existence and Uniqueness For Variational Data Assimilation in Continuous Time}, \href{https://arxiv.org/abs/1805.09269}{ArXiv Preprint 1805.09269}. 

\bibitem{BKO} J. Br\"ocker, T. Kuna, L. Oljaca, \emph{Almost sure error bounds for data assimilation in dissipative systems with unbounded observation noise}, SIAM Journal on Applied Dynamical Systems 17 (4). 2882-2914 (2018).

  
\bibitem{CDP} T. Champion, L. De Pascale, F. Prinari, \emph{$\Ga$-convergence and absolute minimizers for supremal functionals}, COCV ESAIM: Control, Optimisation and Calculus of Variations (2004), Vol. 10, 14-27.

\bibitem{CHY} G. Chen, X. Huang, X. Yang, \emph{Vector Optimization: Set-valued and Variational Analysis}, Lecture Notes in Economics and Mathematical Systems, Springer, 2005.

\bibitem{C} M. G. Crandall, \emph{A visit with the $\infty$-Laplacian}, in \emph{Calculus of Variations and Non-Linear Partial Differential Equations}, Springer Lecture notes in Mathematics 1927, CIME, Cetraro Italy 2005.

\bibitem{CKP} G. Croce, N. Katzourakis, G. Pisante, \emph{$\mathcal{D}$-solutions to the system of vectorial Calculus of Variations in $L^\infty$ via the singular value problem}, Discrete and Continuous Dynamical Systems 37:12, 6165-6181 (2017).

\bibitem{D} B. Dacorogna,  \emph{Direct Methods in the Calculus of Variations}, $2$nd Edition, Volume 78, Applied Mathematical Sciences, Springer, 2008.

\bibitem{HitchhikersSobolevSpaces}
E. Di~Nezza, G. Palatucci, E. Valdinoci,
\newblock Hitchhiker's guide to the fractional Sobolev spaces,
\newblock {\em Bull. Sci. Math.}, 136(5):521--573, 2012.

\bibitem{DPV} M. D'Elia, M. Perego, A. Veneziani, \emph{A Variational Data Assimilation Procedure for the Incompressible Navier-Stokes Equations in Hemodynamics}, Journal of Scientific Computing volume 52, 340-359 (2012).
 
\bibitem{FG} L.C. Florescu, C. Godet-Thobie, \emph{Young measures and compactness in metric spaces}, De Gruyter, 2012.
 
\bibitem{FL} I. Fonseca, G. Leoni, \emph{Modern methods in the Calculus of Variations: $L^p$ spaces}, Springer Monographs in Mathematics, 2007.

\bibitem{FLT} A. Farhat, E. Lunasin, E. S. Titi, \emph{Abridged Continuous Data Assimilation for the 2D Navier-Stokes Equations Utilizing Measurements of Only One Component of the Velocity Field}, J. Math. Fluid Mech. 18, 1-23 (2016).

\bibitem{FMT} C. Foias, C. F. Mondaini, E. S. Titi, \emph{A Discrete Data Assimilation Scheme for the Solutions of the Two-Dimensional Navier-Stokes Equations and Their Statistics}, SIAM J. Appl. Dyn. Syst., 15(4), 2109-2142, 34 pages (2016).

\bibitem{GNP} A. Garroni, V. Nesi, M. Ponsiglione, \emph{Dielectric breakdown: optimal bounds}, Proceedings of the Royal Society A  457, issue 2014 (2001).

\bibitem{Ge} C. Gerhardt, \emph{$L^p$ estimates for solutions to the instationary Navier-Stokes equations in dimension two}, Pacific J. Math. 79:2, 375-398 (1978).

\bibitem{GM} M. Giaquinta, L. Martinazzi, \emph{An Introduction to the Regularity Theory for Elliptic Systems, Harmonic Maps and Minimal Graphs}, Publ. Sc. Norm. Super., vol. 11, Springer, 2012.

\bibitem{Gi} Y. Giga, \emph{Solutions of semilinear parabolic equations in $L^p$ and regularity of weak solutions of the Navier-Stokes equations}, J. Diff. Equations 61, 186-212 (1982).

\bibitem{GS} Y. Giga, H. Sohr, \emph{Abstract $L^p$ estimates for the Cauchy problem with applications to the Navier-Stokes equations in exterior domains}, Journal of functional analysis 102(1), 72-94 (1991).

\bibitem{K0} N. Katzourakis, \emph{An Introduction to Viscosity Solutions for Fully Nonlinear PDE with Applications to Calculus of Variations in $L^\infty$}, Springer Briefs in Mathematics, 2015, DOI 10.1007/978-3-319-12829-0.

\bibitem{K1} N. Katzourakis, \emph{Generalised solutions for fully nonlinear PDE systems and existence-uniqueness theorems}, Journal of Differential Equations 23, 641-686 (2017).

\bibitem{K2} N. Katzourakis, \emph{An $L^\infty$ regularisation strategy to the inverse source identification problem for elliptic equations}, SIAM Journal Math. Analysis, Vol. 51, No. 2, pp. 1349-1370 (2019).

\bibitem{K3} N. Katzourakis, \emph{A minimisation problem in $L^\infty$ with PDE and unilateral constraints}, ESAIM: Control, Optimisation and Calculus of Variations 26, 60 27pp, (2020).

\bibitem{K4} N. Katzourakis, \emph{Inverse optical tomography through PDE-constrained optimisation in $L^\infty$}, SIAM Journal on Control and Optimization, Vol. 57, No. 6, pp. 4205-4233 (2019).

\bibitem{KM} N. Katzourakis, R. Moser, \emph{Existence, Uniqueness and Structure of Second Order Absolute Minimisers}, Archives for Rational Mechanics and Analysis, published online 06/09/2018, DOI: 10.1007/s00205-018-1305-6.

\bibitem{KP1} N. Katzourakis, T. Pryer, \emph{On the numerical approximation of vectorial absolute minimisers}, NoDEA, to appear. 

\bibitem{KP2} N. Katzourakis, T. Pryer, \emph{On the numerical approximation of $p$-Biharmonic and $\infty$-Biharmonic functions}, Numerical Methods for PDE, Numerical Methods in Partial Differential Equations, 1-26 (2018). 

\bibitem{Ko} P. Korn, \emph{Data assimilation for the Navier-Stokes-$\al$ equations}, Physica D 238, 1957--974, (2009).

\bibitem{KZ} C. Kreisbeck, E. Zappale, \emph{Lower semicontinuity and relaxation of nonlocal $L^\infty$-functionals}, Calculus of Variations and PDE 59 (4), 1-36 (2020).

\bibitem{LP} A. Larios, Y. Pei, \emph{Approximate continuous data assimilation of the 2D Navier-Stokes equations via the Voigt-regularization with observable data}, Evolution Equations \& Control Theory 9:3, 733-751, (2019).

\bibitem{M} R. Moser, \emph{Geroch monotonicity and the construction of weak solutions of the inverse mean curvature flow}, Asian J. Math., 19:357-376 (2015).

\bibitem{MS} R. Moser, H. Schwetlick, \emph{Minimizers of a weighted maximum of the Gauss curvature}, Annals of Global Analysis and Geometry, 41 (2), 199 - 207, 2012.

\bibitem{MWZ} Q. Miao, C. Wang, Y. Zhou, \emph{Uniqueness of Absolute Minimizers for $L^\infty$-Functionals Involving Hamiltonians $H(x,p)$}, Archive for Rational Mechanics and Analysis 223 (1), 141-198 (2017).

\bibitem{PZ} F. Prinari, E. Zappale, \emph{A Relaxation Result in the Vectorial Setting and Power Law Approximation for Supremal Functionals}, J Optim. Theory Appl. 186, 412-452 (2020).

\bibitem{RZ} A.N. Ribeiro, E. Zappale, \emph{Existence of minimisers for nonlevel convex functionals}, SIAM J. Control Opt., Vol. 52, No. 5,  (2014) 3341-3370.

\bibitem{SW} A. Schwarz, R. P. Dwight, \emph{Data Assimilation for Navier-Stokes using the Least-Squares Finite-Element Method}, International Journal for Uncertainty Quantification 8(5), 383-403 (2018).

\bibitem{S} V. A. Solonnikov, \emph{Estimates for solution of nonstationary Navier-Stokes equations}, Journal of Soviet Mathematics 8, pages 467-529 (1977).

\bibitem{So} H. Sohr, \emph{The Navier-Stokes Equations}, An Elementary Functional Analytic Approach, Birkha\"user, Springer Basel 2001.

\bibitem{TriebelInterpolationSpaces}
H. Triebel,
\newblock {\em Interpolation theory, function spaces, differential operators},
  volume~18 of {\em North-Holland Mathematical Library},
\newblock North-Holland Publishing Co., Amsterdam-New York, 1978.

\bibitem{Triebel2002}
H. Triebel,
\newblock Function spaces in Lipschitz domains and on Lipschitz manifolds,
  {C}haracteristic functions as pointwise multipliers,
\newblock {\em Rev. Mat. Complut.}, 15(2):475--524, 2002.


\bibitem{Z} E. Zeidler, \emph{Nonlinear Functional Analysis and its Application III: Variational Methods and Optimization}, Springer-Verlag, 1985.



\end{thebibliography}

\end{document}